\input ./style/arxiv-general.cfg
\documentclass[MSNbibl,number,citesort,seceqn,secthm,dvips]{arxbj}
\makeatletter
   \@ifpackageloaded{graphicx}{}{\usepackage{graphicx}}
\makeatother

\volume{22}
\issue{2}
\pubyear{2016}
\firstpage{774}
\lastpage{793}
\doi{10.3150/14-BEJ675}
\docsubty{FLA}

\makeatletter
\newcommand{\rrvert}{\vert}
\newcommand{\llvert}{\vert}
\newcommand{\eqref}[1]{(\ref{#1})}
\newtheorem{theorem}{Theorem}[section]
\newtheorem{corollary}[theorem]{Corollary}
\newremark{example}[theorem]{Example}
\newtheorem{lemma}[theorem]{Lemma}
\newtheorem{proposition}[theorem]{Proposition}
\newremark{remark}[theorem]{Remark}
\makeatother

\begin{document}
\begin{frontmatter}

\title{Simulation of volatility modulated Volterra processes using
hyperbolic stochastic partial differential equations}
\runtitle{Simulation of VMV processes using HSPDEs}

\begin{aug}
\author[A]{\inits{F.E.}\fnms{Fred Espen}~\snm{Benth}\thanksref{A}\ead[label=e1]{fredb@math.uio.no}}
\and
\author[B]{\inits{H.}\fnms{Heidar}~\snm{Eyjolfsson}\corref{}\thanksref{B}\ead[label=e2]{heidar.eyjolfsson@math.uib.no}}
\address[A]{Department of Mathematics, University of
Oslo, P.O. Box 1053 Blindern, N-0316 Oslo, Norway.
\printead{e1}}
\address[B]{Department of Mathematics, University of Bergen, P.O. Box
7803, N-5020 Bergen, Norway.\\
\printead{e2}}
\end{aug}

%
\received{\smonth{4} \syear{2012}}
%
\revised{\smonth{4} \syear{2014}}

\begin{abstract}
We propose a finite difference scheme to simulate solutions to a
certain type of hyperbolic stochastic partial differential equation
(HSPDE). These solutions can in turn estimate so called volatility
modulated Volterra (VMV) processes and L\'evy semistationary (LSS)
processes, which is a class of processes that have been employed to
model turbulence, tumor growth and electricity forward and spot prices.
We will see that our finite difference scheme converges to the solution
of the HSPDE as we take finer and finer partitions for our finite
difference scheme in both time and space. Finally, we demonstrate our
method with an example from the energy finance literature.
\end{abstract}

\end{frontmatter}

\section{Introduction}

This paper is concerned with developing a finite difference scheme to
simulate the so-called mild solution to a particular hyperbolic
stochastic partial differential equation. Our motivation for
considering this scheme is to explore alternative methods to simulate
so-called volatility modulated Volterra (VMV) processes (for
definition, see \eqref{def:GenXVolt}). Volatility modulated Volterra
processes can be simulated by means of numerical integration, but due
to the integrands depending on the time parameter, numerical
integration is cumbersome since at each time step one needs to perform
a complete re-integration. Thus, we propose an alternative method to
simulate these volatility modulated Volterra processes, as the boundary
solution of a hyperbolic stochastic partial differential equation.

We note that a special type of Volatility modulated Volterra processes
are so-called L\'evy semistationary (LSS) processes, which are
processes that are stationary under a stationarity assumption on the
volatility process. These processes have recently been proposed in the
framework of modelling electricity and commodity prices, see
Barndorff-Nielsen, Benth and Veraart \cite{BN1,BN3,BN2}, although they
were initially employed as modelling tools for turbulence and tumor
growth. It has been pointed out that the class of L\'evy semistationary
processes can indeed catch many of the stylised features, such as
spikes and mean-reversion, that have been observed in electricity and
commodity markets. The mean reversion of L\'evy semistationary
processes is in probability, and the high spikes are facilitated by
the volatility process and jumps in the driving L\'evy process. Thus,
it is highly relevant in the energy setting to have an effective
simulating algorithm for derivative pricing purposes.

Employing the finite difference scheme to simulate volatility modulated
Volterra processes as opposed to numerical integration has the
following advantages. As we have already noted we obtain the volatility
modulated Volterra process as the boundary solution of the stochastic
partial differential equation, but in order to obtain a full trajectory
of the boundary with the finite difference scheme we need to
numerically solve the stochastic partial differential equation on a
triangular grid. Therefore when simulating our trajectory, we get the
solution of the stochastic partial differential equation for free on
the triangular grid. In order to simulate a value in a particular point
$(t+\Delta t,x)$ in the grid, we need to know the values at the
previous time step $(t,x)$ as well as at the spatial step above
$(t,x+\Delta x)$. Given initial and boundary conditions we may even
solve the stochastic partial differential equation recursively on a
rectangular grid. We shall show that under certain conditions, the
finite difference scheme converges to the corresponding mild solution.
Moreover, given a stochastic partial differential equation and a
discretization we give a recipe for quantifying the error of the finite
difference scheme in $L^2(\mathbb{P})$.


The rest of the paper is structured as follows. In Section~\ref{LSS},
we start by introducing volatility modulated Volterra processes and
discussing some preliminary results on them which we shall refer to
later in the paper. While in Section~\ref{LSS:HSPDE}, we proceed to
introduce our hyperbolic stochastic partial differential equation, its
mild solution and how we can obtain the volatility modulated Volterra
process as the boundary of the mild solution, under rather general
conditions. Subsequently in Section~\ref{FD}, we introduce the main
contribution of this paper, namely the finite difference scheme for
simulating the hyperbolic stochastic partial differential equation.
Furthermore in that section we discuss convergence results for the
finite difference scheme. Finally in Section~\ref{num}, we present some
numerical examples from the energy literature using our finite
difference scheme, before reaching our concluding remarks in Section~\ref{con}.

\section{Volatility modulated Volterra processes}
\label{LSS}

Throughout this paper, we shall assume that we are working on a given
filtered probability space $(\Omega,\mathcal{F}, \{\mathcal{F}_t\}_{t
\in\mathbb{R}},\mathbb{P})$ which satisfies the usual conditions,
that is, the
probability space is complete, the $\sigma$-algebras $\mathcal{F}_t$
include all the sets in $\mathcal{F}$ of zero probability and the
filtration $\{\mathcal{F}_t\}_{t \in\mathbb{R}}$ is
right-continuous. Note
that the filtration $\{\mathcal{F}_t\}_{t \in\mathbb{R}}$ is indexed
by $\mathbb{R}$.
Following Barndorff-Nielsen \textit{et al.} \cite{BN1,BN3,BN2} we define the
volatility modulated Volterra process (VMV process henceforth) to be a
process of the type
%
\begin{equation}
\label{def:GenXVolt}
X(t) = \mu+ \int_{-\infty}^t p(t,s)a
\bigl(s^-\bigr)\,\mathrm{d}s + \int_{-\infty}^t g(t,s)\sigma
\bigl(s^-\bigr) \,\mathrm{d}L(s),
\end{equation}
%
for $t \in\mathbb{R}$, where $\mu$ is a constant, $\{L(t)\}_{t \in
\mathbb{R}}$ is a
(two-sided) L\'evy process which is adapted to the filtration $\{
\mathcal{F}_t\}_{t \in\mathbb{R}}$, $g$ and $p$ are real-valued deterministic
kernel functions and $\{\sigma(t)\}_{t \in\mathbb{R}}$ and $\{a(t)\}
_{t \in\mathbb{R}
}$ are c\`adl\`ag processes which are adapted to the filtration $\{
\mathcal{F}_t\}_{t \in\mathbb{R}}$. Here the stochastic integral can
be taken
to be defined in the manner developed by Basse-O'Connor \textit{et al.} \cite{BC10}. However although VMV processes can be defined for a rather big
class of L\'evy processes we shall in fact only be concerned with VMV
processes that are driven by square integrable L\'evy processes. Hence,
due to the left limits in the integrand processes, which imply
predictability, the stochastic integration can also be defined in the
sense of Protter \cite{P05}. In addition to the assumption that $\{L(t)\}_{t \in\mathbb{R}}$ is square integrable, we shall also assume
it to
have a zero mean and thus a martingale. Since if $L$ has a drift we may
observe that
\[
L(t) = mt + \bigl(L(t) - mt\bigr),
\]
where $m = \mathbb{E}[L(1)]$ and $\{L(t) - mt\}_{t \in\mathbb{R}}$
is a square
integrable martingale. Thus, we may rewrite~\eqref{def:GenXVolt} as
\[
X(t) = \mu+ \int_{-\infty}^t \bigl(p(t,s)a\bigl(s^-
\bigr) + mg(t,s)\sigma\bigl(s^-\bigr)\bigr)\,\mathrm{d}s + \int_{-\infty}^t
g(t,s)\sigma\bigl(s^-\bigr) \,\mathrm{d}M(s),
\]
where $M(t) = L(t) - mt$ for all $t \in\mathbb{R}$. Moreover, we
shall assume that
%
\begin{equation}
\label{StochProcCond}
\mathbb{E}\bigl[a(t)^2\bigr] \vee\mathbb{E}\bigl[
\sigma(t)^2\bigr] < C
\end{equation}
for some constant $C \geq1 $ and all $t \in\mathbb{R}$. Using this
assumption, we may conclude by Minkowski's integral inequality and It\^{o} isometry that
\begin{eqnarray*}
\mathbb{E}\bigl[X^2(t)\bigr] &\leq & 3 \biggl(\mu^2 +
\mathbb{E} \biggl[ \biggl(\int_{-\infty}^t p(t,s)a
\bigl(s^-\bigr)\,\mathrm{d}s \biggr)^2 \biggr] + \mathbb{E} \biggl[ \biggl(\int
_{-\infty}^t g(t,s)\sigma\bigl(s^-\bigr) \,\mathrm{d}L(s)
\biggr)^2 \biggr] \biggr)
\\
&\leq& 3C \biggl(\mu^2 + \biggl(\int_{-\infty}^t\bigl|p(t,s)\bigr|\,\mathrm{d}s
\biggr)^2 + \int_{-\infty}^tg^2(t,s)\,\mathrm{d}s
\biggr).
\end{eqnarray*}
Thus the VMV process \eqref{def:GenXVolt} is well defined as an element
in $L^2(\mathbb{P})$ if in addition to fulfilling~\eqref{StochProcCond} the
deterministic kernel functions furthermore fulfill
%
\begin{equation}
\label{KernFuncCond}
p(t,\cdot) \in L^1((-\infty,t))\quad\mbox{and}\quad
g(t,\cdot) \in L^2((-\infty,t))
\end{equation}
for all $t \in\mathbb{R}$. In the sequel, we shall always assume that
conditions \eqref{StochProcCond} and \eqref{KernFuncCond} are
fulfilled, which in turn imply that $X(t) \in L^2(\mathbb{P})$ for all
$t \in\mathbb{R}$.

Of particular interest in many applications is the case when $p(t,s) =
p(t-s)$ and $g(t,s) = g(t-s)$, that is, when
%
\begin{equation}
\label{def:XLSS}
X(t) = \mu+ \int_{-\infty}^t p(t-s)a
\bigl(s^-\bigr)\,\mathrm{d}s + \int_{-\infty}^t g(t-s)\sigma
\bigl(s^-\bigr) \,\mathrm{d}L(s).
\end{equation}
Under the additional conditions that the processes $\{a(t)\}_{t \geq
0}$ and $\{\sigma(t)\}_{t \geq0}$ are stationary, the process \eqref
{def:XLSS} is stationary. In particular, we remark that condition
\eqref
{StochProcCond} holds when $a$ and $\sigma$ are stationary. Hence, like
Barndorff-Nielsen \textit{et al.} \cite{BN1,BN3,BN2}, we shall refer to
processes of the type \eqref{def:XLSS} as L\'evy semistationary
processes (or LSS processes). It is furthermore worth noting that in
the case when $p(t,s) = p(t-s)$ and $g(t,s) = g(t-s)$ condition \eqref
{KernFuncCond} is equivalent to $p \in L^1(\mathbb{R}_+)$ and $g \in
L^2(\mathbb{R}_+)$.

Now under the assumption that the stochastic processes $\{a(t)\}_{t \in
\mathbb{R}}$ and $\{\sigma(t)\}_{t \in\mathbb{R}}$ are independent
to each other and
the driving L\'evy process $\{L(t)\}_{t \in\mathbb{R}}$, we have the following
result for LSS processes of the type (\ref{def:GenXVolt}), which is
based on a result in \cite{BEE}.

\begin{proposition}
\label{prop:BN}
Assume that $\{a(t)\}_{t \in\mathbb{R}}$ and $\{\sigma(t)\}_{t \in
\mathbb{R}}$ are
independent to each other and the driving L\'evy process $\{L(t)\}_{t
\in\mathbb{R}}$. Then it holds for processes of the type \eqref{def:GenXVolt} that
%
\begin{equation}
\label{E1}
\mathbb{E}\bigl[X(t)\bigr] = \int_{-\infty}^t
p(t,s)\mathbb{E}\bigl[a\bigl(s^-\bigr)\bigr]\,\mathrm{d}s
\end{equation}
and
%
\begin{equation}
\label{E2}
\mathbb{E}\bigl[X^2(t)\bigr] = \mathbb{E} \biggl[
\biggl(\mu+ \int_{-\infty}^t p(t,s)a\bigl(s^-\bigr)\,\mathrm{d}s
\biggr)^2 \biggr] + \mathbb{E}\bigl[L^2(1)\bigr]\int
_{-\infty}^t g^2(t,s)\mathbb {E}\bigl[
\sigma\bigl(s^-\bigr)^2\bigr] \,\mathrm{d}s.
\end{equation}
In particular if $\{\sigma(t)\}_{t \in\mathbb{R}}$ is stationary
then it holds that
\[
\mathbb{E}\bigl[X^2(t)\bigr] = \mathbb{E} \biggl[ \biggl(
\mu+ \int_{-\infty}^t p(t,s)a\bigl(s^-\bigr)\,\mathrm{d}s
\biggr)^2 \biggr] + \mathbb{E}\bigl[L^2(1)\bigr]\mathbb{E}
\bigl[\sigma^2(0)\bigr]\int_{-\infty}^t
g^2(t,s) \,\mathrm{d}s.
\]
\end{proposition}

\begin{pf}
The characteristic function of the stochastic integral
\[
\widetilde X(t) = \int_{-\infty}^t g(t,s)\sigma
\bigl(s^-\bigr) \,\mathrm{d}L(s)
\]
may be computed by conditioning on the process $\{\sigma(t)\}_{t \in
\mathbb{R}}$:
\[
\varphi_{\widetilde X(t)}(\theta) = \mathbb{E}\bigl[\exp\bigl(\mathrm {i}\theta
\widetilde X(t)\bigr)\bigr] = \mathbb{E} \biggl[\exp \biggl(\int
_{-\infty}^t \psi\bigl(\theta g(t,s)\sigma \bigl(s^-
\bigr)\bigr)\,\mathrm{d}s \biggr) \biggr],
\]
where $\psi$ is the cumulant of $L(1)$, that is, the log-characteristic
function of $L(1)$ (see Proposition~2.6 in \cite{RR}). We observe that
\[
\mathbb{E}\bigl[\widetilde X(t)\bigr]=-\mathrm{i}\psi'(0)\int
_{-\infty
}^tg(t,s)\mathbb{E}\bigl[\sigma\bigl(s^-\bigr)
\bigr]\,\mathrm{d}s=0
\]
since $\mathbb{E}[L(1)]=\psi'(0)=0$ by assumption. Hence, \eqref{E1} follows.
Furthermore, we find
\[
\mathbb{E}\bigl[\widetilde X^2(t)\bigr]=-\psi''(0)
\int_{-\infty
}^tg^2(t,s)\mathbb{E}\bigl[
\sigma\bigl(s^-\bigr)^2\bigr]\,\mathrm{d}s.
\]
So \eqref{E2} follows by independence of the processes $a$, $\sigma$
and $L$.
\end{pf}

In other words, we know everything there is to know about the second
order structure of VMV processes under the assumption that $a$ and
$\sigma$ are independent to each other and the driving L\'evy process.

In Section~\ref{LSS:HSPDE}, we will describe how one can view VMV
processes \eqref{def:GenXVolt} by processes that solve a particular
stochastic partial differential equation, with given initial and
boundary conditions. The solution to the stochastic partial
differential equation can in turn be estimated numerically by a finite
difference method that we will introduce, which will have the same
initial and boundary conditions. For simulating purposes, the initial
condition must be finite and therefore the following lemma will prove useful.

\begin{lemma}
\label{lem:InitCond}
For given VMV processes,
\[
X_1(t) = \int_{-\infty}^t p(t,s)a
\bigl(s^-\bigr)\,\mathrm{d}s + \int_{-\infty}^t g(t,s)\sigma
\bigl(s^-\bigr) \,\mathrm{d}L(s)
\]
and
\[
X_2(t) = \int_{-\infty}^t q(t,s)a
\bigl(s^-\bigr)\,\mathrm{d}s + \int_{-\infty}^t h(t,s)\sigma
\bigl(s^-\bigr) \,\mathrm{d}L(s),
\]
where $\{a(t)\}_{t \in\mathbb{R}}$ and $\{\sigma(t)\}_{t \in\mathbb
{R}}$ satisfy
condition \eqref{StochProcCond} and the deterministic kernel functions
$p,q,g$ and $h$ are square integrable in the sense of \eqref
{KernFuncCond} it holds that
\[
\mathbb{E}\bigl[\bigl|X_1(t) - X_2(t)\bigr|^2\bigr] =
C\bigl(\bigl\Vert p(t,\cdot) - q(t,\cdot )\bigr\Vert_{L^1((-\infty
,t))}^2 + \bigl\Vert g(t,
\cdot)-h(t,\cdot)\bigr\Vert_{L^2((-\infty,t))}^2\bigr),
\]
for a constant $C > 0$ and all $t \in\mathbb{R}$. In particular when
$p(t,s) =
p(t-s)$, $g(t,s) = g(t-s)$, $q(t,s) = q(t-s)$ and $h(t,s) = h(t-s)$ it
holds that
\[
\mathbb{E}\bigl[\bigl|X_1(t) - X_2(t)\bigr|^2\bigr] =
C\bigl(\Vert p-q\Vert_{L^1(\mathbb{R}_+)}^2 + \Vert g-h\Vert_{L^2(\mathbb{R}_+)}^2
\bigr),
\]
for a constant $C > 0$.
\end{lemma}

\begin{pf}
We may apply Proposition~\ref{prop:BN} to obtain
\begin{eqnarray*}
 \mathbb{E} \bigl[\bigl|X_1(t) - X_2(t)\bigr|^2\bigr]
&=& \mathbb{E} \biggl[ \biggl(\int_{-\infty}^t
\bigl(p(t,s) - q(t,s)\bigr)a\bigl(s^-\bigr)\,\mathrm{d}s \biggr)^2 \biggr]\\
&&{}+\mathbb{E}\bigl[L^2(1)\bigr]\int_{-\infty}^t
\bigl(g(t,s) - h(t,s)\bigr)^2\mathbb{E}\bigl[\sigma\bigl(s^-
\bigr)^2\bigr]\,\mathrm{d}s.
\end{eqnarray*}
Moreover it holds by Minkowski's integral inequality that
\[
\mathbb{E} \biggl[ \biggl(\int_{-\infty}^t
\bigl(p(t,s) - q(t,s)\bigr)a\bigl(s^-\bigr)\,\mathrm{d}s \biggr)^2 \biggr] \leq
\mathbb{E}\bigl[a\bigl(s^-\bigr)^2\bigr] \biggl(\int
_{-\infty}^t \bigl|\bigl(p(t,s) - q(t,s)\bigr)\bigr|\,\mathrm{d}s
\biggr)^2.
\]
Now the result follows by \eqref{StochProcCond}.
\end{pf}

Now for a given VMV process \eqref{def:GenXVolt} satisfying conditions
\eqref{StochProcCond} and \eqref{KernFuncCond} we may employ Lemma~\ref
{lem:InitCond} to approximate it with proper stochastic integrals, that
is, integrals over compact intervals. That is, for a given $t \in
\mathbb{R}$,
let $r < t$ be a constant and consider the truncated kernel functions
$\widetilde p(t,s) = 1_{\{s \geq r\}}p(t,s)$ and $\widetilde g(t,s) =
1_{\{s \geq r\}}g(t,s)$. Then due to \eqref{KernFuncCond} it holds that
%
\begin{eqnarray}
&& \bigl\Vert p(t,\cdot) - \widetilde{p}(t,\cdot)\bigr\Vert_{L^1((-\infty,t))}^2
+ \bigl\Vert g(t,\cdot) - \widetilde g(t,\cdot)\bigr\Vert_{L^2((-\infty,t))}^2
\nonumber
\\[-8pt]
\label{InitCond}\\[-8pt]
\nonumber
&& \quad= \int
_{-\infty}^r \bigl(\bigl|p(t,s)\bigr| + g^2(t,s)
\bigr)\,\mathrm{d}s \downarrow0
\end{eqnarray}
as $r \downarrow-\infty$. So by Lemma~\ref{lem:InitCond} we may
approximate the VMV process \eqref{def:GenXVolt} in a fixed point $t
\in\mathbb{R}$ arbitrarily well by a process
%
\begin{equation}
\label{def:FiniteLSS}
X(t) = \mu+ \int_r^t p(t,s)a
\bigl(s^-\bigr)\,\mathrm{d}s + \int_r^t g(t,s)\sigma
\bigl(s^-\bigr) \,\mathrm{d}L(s),
\end{equation}
where $r < t$.






\section{Modelling VMV processes as boundary solutions to HSPDEs}
\label{LSS:HSPDE}

In this section, we will ``raise'' the dimension of our VMV process
(\ref{def:GenXVolt}) to obtain a stochastic process that can be viewed
as a mild solution of a particular hyperbolic stochastic partial
differential equation (HSPDE henceforth). To this end, we need to
define the HSPDE and the concept of a mild solution. For references on
stochastic partial differential equations and mild solutions, we refer
to \cite{CT,PZ}.

For a given $t_0 \in\mathbb{R}$ let us assume that $\{M_t\}_{t \geq
t_0}$ is
a square integrable c\`adl\`ag martingale on a probability space
$(\Omega,\mathcal{F},\mathbb{P})$ with respect to a filtration
$\{\mathcal{F}_t\}_{t \geq t_0}$ that satisfies the usual conditions. Furthermore,
let $\mathcal{P}$ denote the $\sigma$-algebra of predictable sets on
$[t_0,\infty) \times\Omega$, that is, the smallest $\sigma$-algebra of
subsets of $[t_0, \infty) \times\Omega$ containing all sets of the
form $(s,t] \times B$, where $s,t \geq t_0$ and $B \in\mathcal{F}_s$.
Suppose we have a given Hilbert space of univariate real-valued
functions on $\mathbb{R}_+$, denoted by $F$, and predictable (i.e.,
$\mathcal{P}$-measurable) and adapted mappings $\alpha \dvtx  [t_0,\infty) \times
\Omega\to F$ and $\beta: [t_0,\infty) \times\Omega\to F$. Let us
consider the stochastic partial differential equation
%
\begin{equation}
\label{spde:gen}
\mathrm{d}Y(t) = \bigl(AY(t) + \alpha(t)\bigr)\,\mathrm{d}t + \beta(t)\,\mathrm{d}M(t),
\end{equation}
with the initial condition $Y(t_0) = Y_0$, where $Y_0$ is a square
integrable $\mathcal{F}_{t_0}$-measurable random variable with values
in $F$. Here we assume that $A$ is a (potentially unbounded)
infinitesimal generator of a strongly continuous semigroup $\{S(t)\}_{t
\geq0}$ of bounded operators on the Hilbert space $F$. Where by a
strongly continuous semigroup on $F$ we mean that the family of bounded
linear operators $\{S(t)\}_{t \geq0}$ satisfies the following three conditions:
\begin{longlist}[2.]
\item[1.] $S(0) = I$, where $I$ is the identity operator on $F$,
\item[2.] $S(s) \circ S(t) = S(t+s)$ for all $s,t \geq0$,
\item[3.] $\lim_{t \downarrow0}\Vert S(t)f - f\Vert _F = 0$ for all $f \in F$.
\end{longlist}
Note that a family of bounded linear operators $\{S(t)\dvt  t \geq0\}$ on
a Banach space that satisfies the above conditions is called a
$C_0$-\emph{semigroup}. The domain of $A$,
\[
\mathcal{D}(A) = \biggl\{f \in F \dvt  \lim_{t \downarrow0}
\frac{S(t)
f -
f}{t} \mbox{ exists} \biggr\}
\]
will in general be a proper subset of $F$, but it is always a dense
subset in $F$, and its action on $\mathcal{D}(A)$ is given by
%
\begin{equation}
\label{def:A}
Af = \lim_{t \downarrow0} \frac{S(t) f - f}{t},
\end{equation}
for $f \in\mathcal{D}(A)$. For references on operator semigroups see,
for example, \cite{EN}. Since $A$ is in general not a bounded operator,
the notion of a strong solution in the sense of
\[
Y(t) = Y_0 + \int_{t_0}^t \bigl(AY(s)
+ \alpha(s)\bigr)\,\mathrm{d}s + \int_{t_0}^t\beta(s)\,\mathrm{d}M(s)
\]
may not always make sense, since $Y(s)$ might not be in the domain of
$A$. Therefore, the notion of a mild solution to \eqref{spde:gen} has
been introduced in the literature. A \emph{mild solution} to the
equation~\eqref{spde:gen} is a recast of the differential equation
\eqref{spde:gen}:
%
\begin{equation}
\label{mild:sol}
Y(t) = S(t-t_0) Y_0 + \int
_{t_0}^t S(t-s) \alpha(s) \,\mathrm{d}s + \int
_{t_0}^t S(t-s) \beta(s) \,\mathrm{d}M(s).
\end{equation}
In order for the mild solution to be well defined, we need to impose
some conditions on the coefficient functions of \eqref{spde:gen}. The
first integral $\int_{t_0}^t S(t-s)\alpha(s)\,\mathrm{d}s$ is taken to be defined
as a Bochner integral and is thus well defined if the integrand $s
\mapsto S(t-s)\alpha(s)$ is measurable and
%
\begin{equation}
\label{aFin}
\int_{t_0}^t \bigl\Vert S(t-s)
\alpha(s)\bigr\Vert _F\,\mathrm{d}s < \infty.
\end{equation}
Note that the measurability of the integrand from $([t_0,\infty
),\mathcal{B}([t_0,\infty))$ to $(F,\mathcal{B}(F))$ follows from the
strong continuity of the operator semigroup. As for the stochastic
integral $\int_{t_0}^t S(t-s)\beta(s) \,\mathrm{d}M(s)$, recall that by the
Doob--Meyer decomposition, for each c\`adl\`ag square integrable
martingale $\{M(t)\}_{t \geq t_0}$ there exists a unique increasing
predictable process, called the \emph{angle bracket} of $M$, denoted by
$\{\langle M \rangle(t)\}_{t \geq t_0}$ such that $\langle M \rangle
(t_0) = 0$ and $\{M^2(t) - \langle M \rangle(t)\}_{t \geq t_0}$ is a
martingale. For predictable integrands $\beta$ the following It\^o
isometry holds:
\[
\mathbb{E} \biggl[\biggl\Vert \int_{t_0}^t
\beta(s)\,\mathrm{d}M(s)\biggr\Vert _F^2 \biggr] = \mathbb{E} \biggl[\int
_{t_0}^t \bigl\Vert \beta(s)\bigr\Vert _F^2
\,\mathrm{d}\langle M \rangle(s) \biggr],
\]
see, for example, \cite{PZ}. Now for the stochastic integral to be well
defined we need to ensure that the integrand $s \mapsto S(t-s) \beta
(s)$ is predictable, which follows from the strong continuity of the
semigroup, and that it is an element of the space of integrands, that
is, that
%
\begin{equation}
\label{bFin}
\mathbb{E} \biggl[\biggl\Vert \int_{t_0}^t
S(t-s)\beta(s)\,\mathrm{d}M(s)\biggr\Vert_F^2 \biggr] = \mathbb{E} \biggl[
\int_{t_0}^t \bigl\Vert S(t-s)\beta(s)\bigr\Vert _F^2\,\mathrm{d}
\langle M \rangle(s) \biggr] < \infty,
\end{equation}
for all $t \geq t_0$.


In what follows, when we work with solutions to \eqref{spde:gen} we
will always mean mild solutions of the above type (\ref{mild:sol}). Now
let us reconsider the VMV Volterra model (\ref{def:GenXVolt}). Notice
that this equation bears a resemblance to the mild solution (\ref{mild:sol}). However, there are some differences which need to be addressed.

First of all, we need to make an assumption on the operator semigroup
$\{S(t)\}_{t \geq0}$ which is present in \eqref{mild:sol} and the
function space it operates on. Our assumption will be that $\{S(t)\}_{t
\geq0}$ is the strongly continuous semigroup of (left) translation
operators on $F$, defined by
%
\begin{equation}
\label{LShift}
\bigl(S(t)f\bigr) (x) = f(t+x)
\end{equation}
for all $f \in F$ and $x \geq0$. In this case, it follows from \eqref
{def:A} that $A = \partial/\partial x$ is a differential operator on
$F$. Clearly this operator semigroup fulfils the first two algebraic
conditions regardless of the selection of the Hilbert space $F$.
Whereas the third condition by contrast is a topological one, and thus
dependent upon the norm of the Hilbert space. In our setting we propose
to use as state space a Hilbert space proposed by Filipovi\'c \cite
{Fil} in the setting of HJM \cite{HJM} dynamics. For a positive
increasing function $w\dvtx \mathbb{R}_+ \to\mathbb{R}_+$, such that
$\int_0^\infty
w(x)^{-1}\,\mathrm{d}x < \infty$ it is defined as the space of absolutely
continuous functions $f\dvtx  \mathbb{R}_+ \to\mathbb{R}$ satisfying
\[
\int_0^\infty f'(x)^2w(x)\,\mathrm{d}x
< \infty,
\]
endowed with the inner product
\[
\langle f, g \rangle_w = f(0)g(0) + \int_0^\infty
f'(x)g'(x)w(x)\,\mathrm{d}x.
\]
It is easy to see that the norm induced by this inner product satisfies
the strong continuity condition. Moreover for a given $x \geq0$, it
holds that the evaluation functional $\delta_x \dvtx  F \to\mathbb{R}$
defined by
$\delta_x(f) = f(x)$ is uniformly bounded, see \cite{CT}. This will in
turn allow us to evaluate the mild solution \eqref{mild:sol} in any
point $x \geq0$, provided \eqref{aFin} and \eqref{bFin} hold.

The second issue is that the VMV process is defined on $\mathbb{R}$,
whereas a
mild solution to a HSPDE is only defined on a half line $[t_0,\infty)$.
For our purposes, we simply cut the domain of the VMV process in the
following way. For a given $t \in\mathbb{R}$, we assume that there
exists a
$t_0 < t$ such that we can approximate $X(t)$ in \eqref{def:GenXVolt}
by the process
%
\begin{equation}
\label{def:X0}
X_{t_0}(t) = \mu+ \int_{t_0}^t
p(t,s)a\bigl(s^-\bigr)\,\mathrm{d}s + \int_{t_0}^t g(t,s)
\sigma\bigl(s^-\bigr) \,\mathrm{d}L(s)
\end{equation}
in $L^2(\mathbb{P})$. Lemma~\ref{lem:InitCond} confirms that this is possible.
Now having truncated the integration domain, we would like to think of
\eqref{def:X0} as the boundary of a mild solution \eqref{mild:sol}. By
adding a spatial component, $x$, to the above equation we get something
which we may interpret as a mild solution of a HSPDE under assumption
\eqref{LShift} on the operator semigroup. For $x \geq0$, we raise the
dimension of the truncated VMV Volterra model by considering the field
%
\begin{equation}
\label{def:YVolt}
Y(t,x) = \mu+ \int_{t_0}^t
p(t+x,s)a\bigl(s^-\bigr)\,\mathrm{d}s + \int_{t_0}^t g(t+x,s)
\sigma\bigl(s^-\bigr) \,\mathrm{d}L(s).
\end{equation}
Now notice that the process $\{Y(t,\cdot)\}_{t \geq t_0}$ can be viewed
as a mild solution to the HSPDE \eqref{spde:gen}, where $\alpha(t) =
p(t+\cdot,t)a(t-)$, $\beta(t) = g(t+\cdot,t)\sigma(t-)$, $A =
\partial
/\partial x$, $M = L$, and $Y_0 = 0$. Indeed by considering these
coefficient functions for the HSPDE \eqref{spde:gen} under assumption
\eqref{LShift} on the operator semigroup one obtains that the mild
solution \eqref{mild:sol} of the HSPDE \eqref{spde:gen} and the process
defined in \eqref{def:YVolt} coincide. Thus, the VMV Volterra process
\eqref{def:X0} is the boundary solution to the HSPDE, which in turn
approximates the general VMV process \eqref{def:GenXVolt}.

Given the proposed function space selection let us recall the
integrability conditions \eqref{aFin}~and~\eqref{bFin} and inspect what
they translate into in the case of VMV processes. Let us for simplicity
focus on the stochastic integral. In the case, when $M=L$ is a L\'evy
process it holds that $\langle L \rangle(t) = C_1t$, where $C_1 =
\operatorname{Var}[L(1)] > 0$ is a constant. If we furthermore recall the square
integrability condition \eqref{StochProcCond} on the volatility,
condition \eqref{bFin} reduces to
\[
\int_{t_0}^t \bigl\Vert g(t+
\cdot,s)\bigr\Vert_w^2\,\mathrm{d}s = \int_{t_0}^t
g^2(t,s) \,\mathrm{d}s + \int_{t_0}^t\!\! \int
_0^\infty\bigl(g'(t+x,s)
\bigr)^2 w(x) \,\mathrm{d}x\,\mathrm{d}s < \infty.
\]
In particular, this implies that $\Vert g(t,\cdot)\Vert _{L^2((t_0,t))} <
\infty
$ holds. Further strengthening the condition by letting $t_0 \downarrow
-\infty$ and assuming that
%
\begin{equation}
\label{KernFuncSup}
\int_{-\infty}^t\bigl\Vert g(t+
\cdot,s)\bigr\Vert _w^2\,\mathrm{d}s < \infty
\end{equation}
implies that the condition \eqref{KernFuncCond} is satisfied by $g$.
Thus we observe that assuming that \eqref{KernFuncSup} holds for $g$
and $p$, is sufficient for our purposes with VMV processes and for the
solution of the HSPDE to be well defined.

In what follows, for a given discretization $t_1 < t_2 < \cdots< t_N$
in the time domain, to simulate a trajectory $\{X(t_n)\}_{n=0}^N$ of
the VMV Volterra process \eqref{def:GenXVolt}, we propose the following
two step procedure:
\begin{longlist}[2.]
\item[1.]
Truncate the integration domain of \eqref{def:GenXVolt} from $\mathbb
{R}$ to
$[t_0, \infty)$. Where $t_0 \leq t_1$ is such that $\Vert X(t) -
X_{t_0}(t)\Vert _{L^2(\mathbb{P})}$ is close to zero for $t \geq t_1$.

\item[2.] Raise the dimension of the truncated VMV Volterra model by
considering the field \eqref{def:YVolt}. Now simulate the HSPDE \eqref
{spde:gen} with $\alpha(t) = p(t+\cdot,t)a(t-)$, $\beta(t) =
g(t+\cdot
,t)\sigma(t-)$, $M = L$ and $Y_0=0$ under assumption \eqref{LShift} on
the operator semigroup using the finite difference scheme that will be
introduced in Section~\ref{FD}. The trajectory $\{X(t_n)\}_{n=0}^N = \{
Y(t_n,0)\}_{n=0}^N$ is obtained as the boundary solution of the HSPDE.
\end{longlist}


In many cases, one may even be interested in more than just the
boundary, as the following example shows.

\begin{example}
\label{ex:forward}
In Section~5 of \cite{BN2} Barndorff-Nielsen \textit{et al.}, derive a model for
pricing electricity forward contracts based on general L\'evy driven
Volterra electricity spot prices. Thus deseasonalized electricity spot
prices $\{X(t)\}_{t \in\mathbb{R}}$ are generally modelled by VMV
processes of
the type
\[
X(t) = \int_{-\infty}^t g(t,s)\sigma\bigl(s^-\bigr)
\,\mathrm{d}L(s),
\]
where the components of the integral fulfil all the necessary
conditions listed in Section~\ref{LSS}. Examples of kernel functions
considered by Barndorff-Nielsen \textit{et al.}~\cite{BN2} include $g(t,s) =
\exp
(-\alpha(t-s))$ for a constant $\alpha> 0$, and $g(t,s) = \sigma
/(t-s+b)$ for constants $\sigma, b > 0$. Under certain integrability
conditions, forward price dynamics $F_t(T)$ may be derived as an
expression involving the volatility modulated Volterra process
\[
\int_{-\infty}^t g(T,s)\sigma\bigl(s^-\bigr) \,\mathrm{d}L(s).
\]
Here $T$ is time of delivery. Letting $x = T-t$ we may write
\[
\int_{-\infty}^t g(T,s)\sigma\bigl(s^-\bigr) \,\mathrm{d}L(s) =
\int_{-\infty}^t g(t+x,s)\sigma\bigl(s^-\bigr) \,\mathrm{d}L(s),
\]
and we are back to our mild solution. Hence, in a practical context, we
are interested in simulating the joint spot-forward price dynamics.
This can be done by simulating the mild solution of the corresponding
HSPDE. Hence, in an energy market context, the finite difference scheme
approach gives a joint simulation of spot and forward prices for all
maturities directly without re-integration at each maturity.
\end{example}

We shall return to this example in Section~\ref{num}, after we have
discussed our finite difference scheme.



\section{The finite difference scheme}
\label{FD}
This section presents the main contribution of this paper, namely a
finite difference scheme for simulating solution fields for the HSPDE
\eqref{spde:gen}, under the assumption $A = \partial/\partial x$.


Now let us introduce the following notation for the finite difference
method. Let $\Delta x > 0$ and $\Delta t > 0$ denote the discrete steps
in space and time respectively, and denote by
\[
y_j^n \approx Y(t_0+n\Delta t) (j\Delta x)
\]
the approximation of the solution of \eqref{spde:gen} at the point
$(t_0 + n\Delta t, j\Delta x)$, where $n=0,\ldots,N$ and $j=0,\ldots,J$
for some $J,N \in\mathbb{N}$. From our HSPDE \eqref{spde:gen} with
$A=\partial/\partial x$, using forward finite difference, that is, by using the
approximations $\mathrm{d}Y(t) \approx Y(t+\Delta t) - Y(t)$, $\mathrm{d}t \approx\Delta
t$, $\mathrm{d}M(t) \approx M(t+\Delta t) - M(t)$ and $AY(t) \approx(Y(t)(\cdot
+ \Delta x) - Y(t))/\Delta x$, we derive the finite difference scheme
%
\begin{equation}
\label{def:FD} y_j^{n+1} = \lambda y_{j+1}^n
+ (1-\lambda)y_j^n + \alpha_j^n
\Delta t + \beta_j^n \Delta M^n,
\end{equation}
where $\lambda= \Delta t / \Delta x$, $x_j = j\Delta x$, $t_n = t_0 +
n\Delta t$, $\alpha_j^n = \alpha(t_n)(x_j)$, $\beta_j^n = \beta
(t_n)(x_j)$ and $\Delta M^n = M(t_{n+1}) - M(t_n)$. Clearly, one should
adjust the initial value so that it fits with the initial value of the
HSPDE one is interested in simulating, that is, by setting $y_j^0 =
Y_0(x_j)$ for all $j=0,\ldots,J$. For instance in our VMV applications
(recall \eqref{def:YVolt}) this means letting $y_j^0 = \mu$, for all
$j=0,\ldots,J$. We furthermore note that information about the initial
values are sufficient. Since in order to obtain a value at a given
point $(t_{n+1},x_j)$ the scheme requires information about the values
at the previous time steps $(t_n,x_j)$ and $(t_n,x_{j+1})$. Thus for a
fixed $j'$ in order to calculate the trajectory $\{y_{j'}^n\}_{n=1}^N$
we only need information about the previous values on a triangular
grid, that is, we need to know the values of
%
\begin{equation}
\label{triang} %
\begin{array} {l} y_{j'+N}^0;
\\[6pt]
y_{j'+N-1}^0,y_{j'+N-1}^{1};
\\[2pt]
\vdots
\\[2pt]
y_{j'+1}^0, y_{j'+1}^1,
\ldots,y_{j'+1}^{N-1},
\end{array} %
\end{equation}
all of which may be obtained from the initial values. Hence, to
simulate the random field $\{Y(t_n)(x_j)\}_{j=0,n=0}^{J,N}$ which is
the solution of the HSPDE \eqref{spde:gen} on a rectangular grid, for a
given initial value, without knowing the values at the boundary
$(x_J)$, using the finite difference scheme~\eqref{def:FD}, we propose
the following.
\begin{longlist}[3.]
\item[1.] Simulate $\Delta M^n$, for $n=0,\ldots,N-1$.
\item[2.] Compute the values of the triangular grid \eqref{triang} where $j'=J$.
\item[3.] Compute the values of the rectangular grid, using values from the
triangular grid where necessary.
\end{longlist}
%
We remark that in some cases it may however be natural to impose a
boundary condition on the spatial dimension. In the case of LSS
processes with $p(t,s) = p(t-s)$ and $g(t,s) = g(t-s)$ one could use
Lemma~\ref{lem:InitCond} to assume $y_J^n = 0$, for all $n=0,\ldots,N$,
if $x_J$ is big enough, since \eqref{KernFuncCond} implies that $p \in
L^1(\mathbb{R}_+)$ and $g \in L^2(\mathbb{R}_+)$, so they vanish at infinity.

As in the case of a finite difference scheme for the standard advection
partial differential equation, one needs some constraints on the
discrete steps, that is, $(\Delta x,\Delta t)$, to guarantee its
stability. The stability condition of Courant, Friedrichs, and Lewy
(the CFL condition, see \cite{CFL}) is needed to ensure the stability
of our finite difference scheme \eqref{def:FD}. In our case this
translates into the necessary constraint
%
\begin{equation}
\label{stable}
\Delta t \leq\Delta x,
\end{equation}
which we assume to hold.

For the rest of this section, we will study the convergence properties
of the finite difference scheme. Given our function space $F$ of
real-valued functions equipped with a supremum norm it will be
convenient for our analysis to define the following family of bounded
linear operators on~$F$. Given positive $\Delta x > 0$ and $\Delta t >
0$ corresponding to the steps of the finite difference scheme in space
and time respectively let us consider the family $\{T_{\Delta x,\Delta
t}\}_{\Delta x > 0, \Delta t > 0}$ which is defined by
%
\begin{equation}
\label{def:T} T_{\Delta x,\Delta t} = I + \Delta t\frac{S(\Delta x) - I}{\Delta x},
\end{equation}
for all $\Delta x > 0, \Delta t > 0$, where $I$ denotes the identity
operator on $F$ and $S(\Delta x)$ is the left shift operator whose
action on $F$ is given by \eqref{LShift}. The following lemma will be
useful for proving convergence of the finite difference scheme.

\begin{lemma}
\label{lem:itscheme}
For given steps $\Delta x > 0$ in space and $\Delta t > 0$ in time, the
finite difference scheme~\eqref{def:FD} admits the representation
%
\begin{equation}
\label{y:it} y^n_j = T^n
y^0_j + \sum_{i=0}^{n-1}
T^{n-1-i} \alpha^i_j \Delta t + \sum
_{i=0}^{n-1} T^{n-1-i} \beta^i_j
\Delta M^i
\end{equation}
for all $n=0,\ldots,N$ and $j=0,\ldots,J$, where $T = T_{\Delta x,
\Delta t}$ is defined by \eqref{def:T} and where $T^n = T^{\circ n}$
denotes the composition of the operator $T$ with itself $n$ times and
$T^0 = I$.
\end{lemma}

\begin{pf}
We proceed by means of induction on $n$. The identity (\ref{y:it})
clearly holds for $n=0$ and all $j=0,\ldots,J$. Supposing that the
identity (\ref{y:it}) is satisfied by some $n \geq0$ and all
$j=0,\ldots,J$, we obtain the following.
\begin{eqnarray*}
y^{n+1}_j &= &\lambda y^n_{j+1} + (1-
\lambda)y^n_j + \alpha^n_j
\Delta t + \beta^n_j \Delta M^n
\\
&= & T^ny^0_j + \lambda\bigl(T^ny^0_{j+1}
- T^ny^0_j\bigr)
\\
&&{}+ \sum_{i=0}^{n-1} \bigl(T^{n-1-i}
\alpha^i_j + \lambda \bigl(T^{n-1-i}
\alpha^i_{j+1} - T^{n-1-i}\alpha^i_j
\bigr) \bigr)\Delta t + \alpha ^n_j\Delta t
\\
&&{}+ \sum_{i=0}^{n-1} \bigl(T^{n-1-i}
\beta^i_j + \lambda \bigl(T^{n-1-i}\beta
^i_{j+1} - T^{n-1-i}\beta^i_j
\bigr) \bigr)\Delta M^i + \beta^n_j\Delta
M^n
\\
&=& TT^n y^0_j + \sum
_{i=0}^{n-1}\bigl(TT^{n-1-i}
\alpha^i_j\bigr)\Delta t + \alpha ^n_j
\Delta t + \sum_{i=0}^{n-1}
\bigl(TT^{n-1-i}\beta^i_j\bigr)\Delta
M^i + \beta ^n_j \Delta M^n
\\
&= & T^{n+1} y^0_j + \sum
_{i=0}^n T^{n-i} \alpha^i_j
\Delta t + \sum_{i=0}^n T^{n-i}
\beta^i_j \Delta M^i.
\end{eqnarray*}
This completes the proof.
\end{pf}

The above lemma characterizes the finite difference scheme (\ref
{def:FD}) for a given discretization as the sum of three entities
which, under appropriate conditions, will converge to their
corresponding parts in the mild solution (\ref{mild:sol}) as we
consider finer and finer partitions in time and space. More precisely,
we will employ the fact that the composed operator $T^n$, where
$T=T_{\Delta x, \Delta t}$ is defined by (\ref{def:T}), converges to
the left shift operator $S(t_n-t_0)$ as we consider finer and finer
partitions in first time and then space.

Let us take a closer look on the family \eqref{def:T} of operators. The
following lemma will be employed later for proving a convergence result
on the finite difference scheme.

\begin{lemma}
\label{lem:Tconv}
Suppose $\beta\dvtx  \Omega\times\mathbb{R}_+ \to\mathbb{R}$ is a
function that satisfies
the Lipschitz condition
\[
\mathbb{E}\bigl[\bigl|\beta(x_1) - \beta(x_2)\bigr|^2
\bigr] \leq L|x_1 -x_2|^2
\]
for all $x_1,x_2 \geq0$ where $L > 0$ is a constant. Then
\[
\mathbb{E}\bigl[\bigl|T^m \beta(x) - S(t)\beta(x)\bigr|^2\bigr]
\leq Lt(\Delta x - \Delta t),
\]
where $\{S(t)\}_{t \geq0}$ denotes the left shift semigroup \eqref
{LShift} and $T$ is defined in \eqref{def:T} with $\Delta t = t/m$ and
$\Delta t \leq\Delta x$, for all $x \geq0$, $t > 0$ and $m \geq1$.
\end{lemma}

\begin{pf}
Let $\lambda= \Delta t/\Delta x$ and suppose first that $\lambda= 1$,
then clearly $T = S(\Delta x)$ and $T^m = S(t)$. Now suppose that
$\lambda< 1$, and observe that by the binomial theorem it holds that
\begin{eqnarray*}
T^m\beta(x) &=& (1-\lambda)^m \biggl(I +
\frac{\lambda}{1-\lambda
}S(\Delta x) \biggr)^m\beta(x)
\\
&=&  \sum_{k=0}^m {m \choose k}
\lambda^k (1-\lambda)^{m-k}\beta(x + k\Delta x) =
\mathbb{E}'\bigl[\beta(x+\Delta xZ)\bigr],
\end{eqnarray*}
where $Z$ denotes a binomial random variable with parameters $m$ and
$\lambda$ on a probability space $(\Omega',\mathcal{F}',\mathbb
{P}')$, with
expectation operator denoted by $\mathbb{E}'$. Now recall that a binomial
random variable has expected value $m\lambda$ and variance $m\lambda
(1-\lambda)$, from which it is easy to deduce that the random variable
$\Delta xZ$ has expected value $t$ and variance $t(\Delta x -\Delta
t)$. Thus by the Cauchy--Schwarz inequality, the Fubini theorem (or the
linearity of the expected value) and the Lipschitz condition
\begin{eqnarray*}
\mathbb{E}\bigl[\bigl|T^m\beta(x) - \beta(x+t)\bigr|^2\bigr] &
\leq & \mathbb{E}\bigl[\mathbb {E}'\bigl[\bigl|\beta(x + \Delta xZ) -
\beta(x+t)\bigr|^2\bigr]\bigr]
\\
&\leq&  L\mathbb{E}'\bigl[|\Delta xZ - t|^2\bigr]
\\
&=& Lt(\Delta x - \Delta t).
\end{eqnarray*}
This concludes the proof.
\end{pf}

In the context of the finite difference scheme \eqref{def:FD}, $t-t_0$
corresponds to the length of the time interval $[t_0,t]$ and taking
$\Delta t = t_{n+1} - t_n$ for all $n=0,\ldots,N-1$. In light of Lemma~\ref{lem:Tconv}, it is interesting to comment on the difference between
employing our finite difference scheme as opposed to numerical
integration. Supposing that we are mainly interested in the boundary
solution ($x=0$) of a particular HSPDE on a given time grid $t_0 <
\cdots< t_N = t$, then one could estimate its mild solution at a
particular time step $t$ by means of numerical integration in the
following way:
%
\begin{equation}
\label{NumInt} 
\hspace*{-10pt}\widetilde Y(t) (0) = S(t-t_0)Y_0(0) +
\sum_{i=0}^{N-1} S(t-t_{i+1})
\alpha(t_i) (0)\Delta t + \sum_{i=0}^{N-1}
S(t-t_{i+1})\beta (t_i) (0)\Delta M^i.
\end{equation}
By comparison to \eqref{y:it}, one sees that the above equation is
quite similar. Moreover, Lemma~\ref{lem:Tconv} provides us with some
evidence that the equation \eqref{y:it} for $j=0$ and the above
equation \eqref{NumInt} tend to give us the same trajectories as we
consider finer and finer steps. In particular when $\Delta t = \Delta
x$ the two approaches give us the exact same trajectories. But the
difference between the two respective methods, given that the
coefficient functions are sufficiently well behaved, is that one of
them only gives us the boundary solution of the HSPDE, whereas the
other one solves the HSPDE on a triangular grid.

This is relevant in the setting of Example~\ref{ex:forward}, in the
context of simulating the joint spot-forward dynamics. Another
advantage of employing the finite difference, is that given the values
$Y(t)(0)$ and $Y(t)(\Delta x)$ for a particular $t > t_0$ we easily
obtain the next value $Y(t+\Delta t)(0)$ by means of the finite
difference scheme \eqref{def:FD}. However if we employ numerical
integration we cannot use this information to calculate the next step
$Y(t+\Delta t)(0)$, we need to do a complete re-integration in the time domain.

We have the following convergence result, which can be used to
determine whether or not the finite difference scheme \eqref{def:FD} is
convergent in $L^2(\mathbb{P})$ for a particular HSPDE and to
determine the
convergence rate. We shall only consider HSPDEs with initial value and
coefficient functions that are uniformly Lipschitz in the sense of
Lemma~\ref{lem:Tconv}. That it we assume that
%
\begin{eqnarray}
&& \mathbb{E}\bigl[\bigl|Y_0(x_1) -
Y_0(x_2)\bigr|^2\bigr] \vee\mathbb{E}\bigl[\bigl|
\alpha(s) (x_1) - \alpha (s) (x_2)\bigr|^2\bigr]
\vee\mathbb{E}\bigl[\bigl|\beta(s) (x_1) - \beta(s) (x_2)\bigr|^2
\bigr]\hspace*{8pt}
\nonumber
\\[-8pt]
\label{coeffLip}\\[-8pt]
\nonumber
&& \quad \leq L|x_1 -x_2|^2
\end{eqnarray}
hold for all $s \in[t_0,t]$, $x_1,x_2 \geq0$ and a constant $L > 0$.


\begin{proposition}
\label{prop:Conv}
Consider the finite difference scheme \eqref{def:FD} under the
representation \eqref{y:it}, where the initial value and the
coefficient functions satisfy the Lipschitz condition \eqref{coeffLip}.
Suppose furthermore that the coefficient functions are independent of
the driving martingale process. Then if $t_n = t_0 + n\Delta t$ and
$x_j = j\Delta x$, for $n,j \geq0$, it holds that
\begin{eqnarray*}
\mathbb{E} \bigl[\bigl\llvert y_j^n -
Y(t_n) (x_j)\bigr\rrvert ^2 \bigr]  &\leq &
C_1(n) (\Delta x - \Delta t) + C_2(n)\Delta
t^2
\\
&& {}+ C_3(n)\mathbb{E} \Bigl[\sup_{0 \leq s-r < \Delta
t}\bigl|S(t_n-s)
\bigl(\alpha(r) - \alpha(s)\bigr)\bigr|^2 \Bigr]
\\
&&{}+ C_4(n)\mathbb{E} \Bigl[\sup_{0 \leq s-r < \Delta
t}\bigl|S(t_n-s)
\bigl(\beta(r) - \beta(s)\bigr)\bigr|^2 \Bigr],
\end{eqnarray*}
where
\begin{eqnarray*}
C_1(n) &=&  3L(t_n - t_0) \bigl\{1 +
4(t_n-t_0)^2 + 4\mathbb{E}\bigl[\langle M
\rangle (t_n)\bigr] \bigr\},\\
 C_2(n) &=&  12L\bigl
\{(t_n-t_0)^2 + \mathbb{E}\bigl[\langle M
\rangle (t_n)\bigr]\bigr\},
\\
C_3(n) &=&  12L(t_n-t_0)^2\quad \mbox{and} \quad C_4(n) = 12\mathbb{E}\bigl[\langle M \rangle(t_n)
\bigr].
\end{eqnarray*}
\end{proposition}

\begin{pf}
First, notice that
\[
\mathbb{E}\bigl[\bigl|T^{N-1-i}Y_0 - S(t-t_0)Y_0\bigr|^2
\bigr] \leq L(t-t_0) (\Delta x - \Delta t)
\]
follows directly from Lemma~\ref{lem:Tconv}. Since $M$ is square
integrable and independent to $\beta$ it holds by It\^o isometry and
Lemma~\ref{lem:Tconv} that
\begin{eqnarray*}
&& \mathbb{E} \Biggl[\Biggl\llvert \sum_{i=0}^{N-1}T^{N-1-i}
\beta(t_i)\Delta M^i - \sum
_{i=0}^{N-1}S(t-t_{i+1})\beta(t_i)
\Delta M^i\Biggr\rrvert ^2 \Biggr]
\\
&&\quad= \mathbb{E} \Biggl[\int_{t_0}^t
\sum_{i=0}^{N-1} \bigl(T^{N-1-i}\beta
(t_i) - S(t-t_{i+1})\beta(t_i)
\bigr)^2 1_{[t_i,t_{i+1})}(s) \,\mathrm{d}\langle M \rangle (s) \Biggr]
\\
&&\quad= \sum_{i=0}^{N-1} \mathbb{E}\bigl[
\bigl(T^{N-1-i}\beta(t_i) - S(t-t_{i+1})\beta
(t_i)\bigr)^2\bigr]\mathbb{E} \biggl[\int
_{t_i}^{t_{i+1}}\,\mathrm{d}\langle M \rangle (s) \biggr]
\\
&&\quad\leq L(\Delta x - \Delta t)\sum_{i=0}^{N-1}
(t-t_{i+1})\mathbb {E}\bigl[\langle M \rangle(t_{i+1}) -
\langle M \rangle(t_i)\bigr] \leq L(t-t_0)\mathbb{E}
\bigl[\langle M \rangle(t)\bigr](\Delta x - \Delta t).
\end{eqnarray*}
Furthermore by Lipschitz continuity and independence of $M$ and $\beta$
we get that
\begin{eqnarray*}
&&\!\! \mathbb{E} \Biggl[\Biggl\llvert \sum_{i=0}^{N-1}S(t-t_{i+1})
\beta (t_i)\Delta M^i - \int_{t_0}^t
S(t-s)\beta(s)\,\mathrm{d}M(s)\Biggr\rrvert ^2 \Biggr]
\\
&&\!\!\quad= \sum_{i=0}^{N-1}
\mathbb{E} \biggl[\int_{t_i}^{t_{i+1}}
\bigl(S(t-t_{i+1})\beta (t_i) - S(t-s)\beta(s)
\bigr)^2\,\mathrm{d}\langle M \rangle(s) \biggr]
\\
&&\!\! \quad\leq\sum_{i=0}^{N-1}\mathbb{E} \Bigl[\sup
_{s \in
[t_i,t_{i+1})} \bigl(S(t-t_{i+1})\beta(t_i) -
S(t-s)\beta(s) \bigr)^2 \Bigr]\mathbb {E} \biggl[\int
_{t_i}^{t_{i+1}}\,\mathrm{d}\langle M \rangle(s) \biggr]
\\
&&\!\!\quad\leq2\sum_{i=0}^{N-1}\mathbb{E} \Bigl[\sup
_{s \in
[t_i,t_{i+1})}\bigl(\bigl(S(t-t_{i+1})\beta(t_i) -
S(t-s)\beta(t_i)\bigr)^2
+ \bigl(S(t-s)\beta(t_i) - S(t-s)\beta(s)
\bigr)^2\bigr) \Bigr]\\
&&\!\!\hspace*{22pt}\qquad {}\times\mathbb{E} \biggl[\int_{t_i}^{t_{i+1}}\,\mathrm{d}
\langle M \rangle(s) \biggr]
\\
&&\!\!\quad\leq 2\mathbb{E}\bigl[\langle M \rangle(t)\bigr] \Bigl(L\Delta t^2 +
\mathbb {E} \Bigl[\sup_{
0\leq s-r < \Delta t}\bigl|S(t-s) \bigl(\beta(r) - \beta(s)
\bigr)\bigr|^2 \Bigr] \Bigr). 
\end{eqnarray*}
Putting the above inequalities together and employing the elementary
inequality $(x+y)^2 \leq2(x^2+y^2)$, we obtain
\begin{eqnarray*}
&& \mathbb{E} \Biggl[\Biggl\llvert \sum_{i=0}^{N-1}T^{N-1-i}
\beta(t_i)\Delta M^i - \int_{t_0}^t
S(t-s)\beta(s)\,\mathrm{d}M(s)\Biggr\rrvert ^2 \Biggr]
\\
&&\quad\leq 4\mathbb{E}\bigl[\langle M
\rangle(t)\bigr] \Bigl(L(t-t_0) (\Delta x - \Delta t) + L\Delta
t^2 + \mathbb{E} \Bigl[\sup_{0 \leq s-r < \Delta
t}\bigl|S(t-s) \bigl(
\beta(r) - \beta(s)\bigr)\bigr|^2 \Bigr] \Bigr).
\end{eqnarray*}
Now
\begin{eqnarray*}
&&\mathbb{E} \Biggl[\Biggl\llvert \sum_{i=0}^{N-1}T^{N-1-i}
\alpha(t_i)\Delta t - \int_{t_0}^t
S(t-s)\alpha(s)\,\mathrm{d}s\Biggr\rrvert ^2 \Biggr]
\\
&&\quad\leq 4(t-t_0)^2 \Bigl(L(t-t_0) (\Delta x -
\Delta t) + L\Delta t^2 + \mathbb{E} \Bigl[\sup_{0 \leq s-r < \Delta t}\bigl|S(t-s)
\bigl(\alpha(r) - \alpha (s)\bigr)\bigr|^2 \Bigr] \Bigr)
\end{eqnarray*}
follows in a similar manner, replacing the It\^o isometry argument by
the Cauchy--Schwarz inequality. The proof is completed by employing the
representation in Lemma~\ref{lem:itscheme}, the elementary inequality
$(x+y+z)^2 \leq3(x^2 + y^2 + z^2)$, and collecting the resulting terms.
\end{pf}


Now let us finish the section by coming back to LSS processes. In this
case, we are only concerned with HSPDEs that have coefficient functions
which can be separated into a stochastic part and a deterministic part.
That is, HSPDEs that have coefficient functions on the following form:
%
\begin{equation}
\label{sepvar}
\alpha(t) = pa(t-)\quad \mbox{and}\quad \beta(t) = g\sigma(t-),
\end{equation}
where $p,g \in F$ are Lipschitz continuous functions with a joint
Lipschitz constant $L > 0$, and $\{a(t)\}_{t \geq t_0}$ and $\{\sigma
(t)\}_{t \geq t_0}$ are predictable and adapted stochastic processes
that satisfy \eqref{StochProcCond}. We shall moreover require that
%
\begin{equation}
\label{BoundedK}
|g|^2 \vee|p|^2 < K
\end{equation}
for a constant $K \geq1$. Indeed for our function space equipped with
a supremum norm these assumptions guarantee that the corresponding
HSPDE has a well-defined mild solution.

\begin{corollary}
Consider the finite difference scheme \eqref{def:FD} under the
representation \eqref{y:it}, where the initial value and the
coefficient functions satisfy the Lipschitz condition \eqref{coeffLip}.
Suppose furthermore that the coefficient functions are independent of
the driving martingale process, and that \eqref{sepvar} and \eqref
{BoundedK} hold. Then if $t_n = t_0 + n\Delta t$ and $x_j = j\Delta x$,
for $n,j \geq0$, it holds that
\begin{eqnarray*}
&& \mathbb{E} \bigl[\bigl\llvert y_j^n -
Y(t_n) (x_j)\bigr\rrvert ^2 \bigr] \\
&&\quad\leq
C_1(n) (\Delta x - \Delta t) + C_2(n)\Delta
t^2
\\
&&\qquad{}+ C_3(n)\mathbb{E} \Bigl[\sup_{0 \leq s-r < \Delta t}\bigl|a(r) -
a(s)\bigr|^2 \Bigr] + C_4(n)\mathbb{E} \Bigl[\sup
_{0 \leq s-r < \Delta t}\bigl|\sigma(r) - \sigma (s)\bigr|^2 \Bigr],
\end{eqnarray*}
where
\begin{eqnarray*}
C_1(n) &=&  3L(t_n - t_0) \bigl\{1 +
4(t_n-t_0)^2 + 4\mathbb{E}\bigl[\langle M
\rangle (t_n)\bigr] \bigr\}, \\
 C_2(n) &=&  12L\bigl
\{(t_n-t_0)^2 + \mathbb{E}\bigl[\langle M
\rangle (t_n)\bigr]\bigr\},
\\
C_3(n) &=&  12KL(t_n-t_0)^2\quad
\mbox{and}\quad C_4(n) = 12K\mathbb {E}\bigl[\langle M
\rangle(t_n)\bigr].
\end{eqnarray*}
In particular if $M$ is a L\'evy process then $\mathbb{E}[\langle M
\rangle
(t)] = Ct$ for a constant $C \geq0$.
\end{corollary}

\section{Numerical examples}
\label{num}

In this section, we present some numerical examples to illustrate the
finite difference scheme and our convergence results in the previous
section. As an example, consider
%
\begin{equation}
\label{modBj}
g(u) = \frac{a}{u+b}\mathrm{e}^{-\alpha u},
\end{equation}
where $a,b > 0$ and $\alpha\geq0$. This is a blend of the kernel
function suggested by Bjerksund \textit{et al.} \cite{BRS} and the OU process,
and thus constitutes a potential kernel function for applications in
electricity. Returning to Example~\ref{ex:forward}, for a fixed grid in
time $t_0 < t_1 < \cdots< t_N$ and space $0=x_0<x_1 < \cdots< x_J$
with fixed increments $\Delta t$ and $\Delta x$ respectively, consider
simulating the random field
%
\begin{equation}
\label{GenBjerk}
Y(t,x) = \int_0^t g(t-s+x)
\sigma\bigl(s^-\bigr)\,\mathrm{d}B(s),
\end{equation}
where $g$ represents the kernel function \eqref{modBj}, $B$ is standard
Brownian motion and $\sigma^2(t) = Z(t)$, where
%
\begin{equation}
\label{sigma2} Z(t) = \int_{-\infty}^t
\mathrm{e}^{-\lambda(t-s)}\,\mathrm{d}U(s),
\end{equation}
and $U$ is a subordinator process.
%
\begin{figure}

\includegraphics{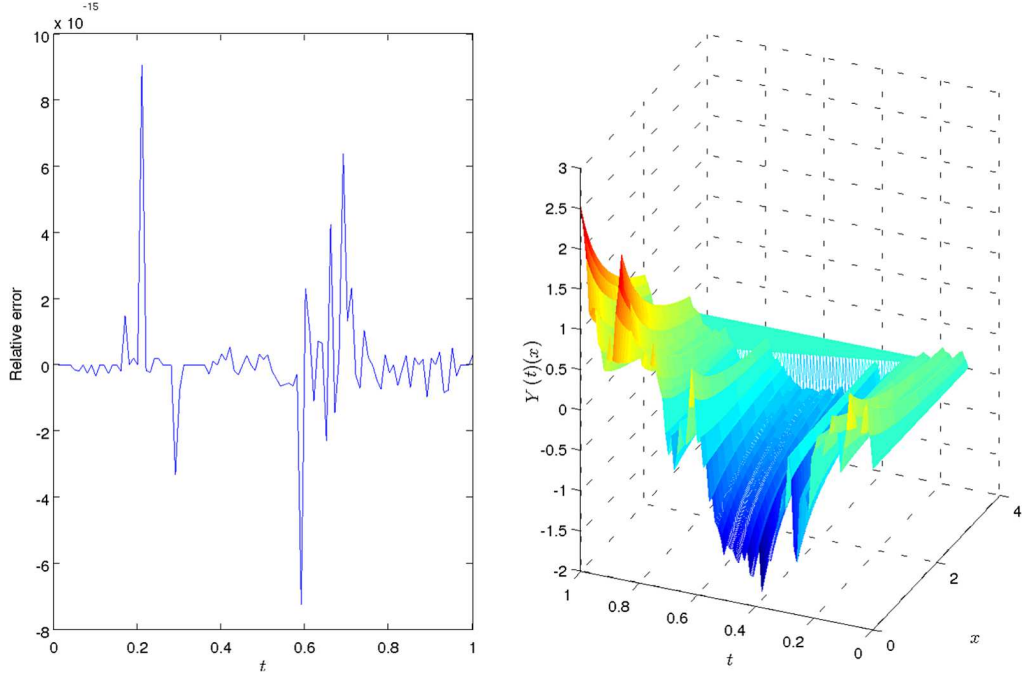}

\caption{Left: the relative error of the boundary path of \protect\eqref{GenBjerk} ($x=0$) where $g$ is given by \protect\eqref{modBj}, obtained by
numerical integration versus the finite difference scheme. Right: the
field \protect\eqref{GenBjerk} where $g$ is given by \protect\eqref{modBj}, obtained by
the finite difference method, on a rectangular grid with step sizes
$\Delta t = \Delta x = 0.01$.}
\label{figIGGB}
\end{figure}
Now simulating \eqref{GenBjerk} on a rectangular grid with the finite
difference method is much more efficient than using numerical
integration to calculate each trajectory for a fixed $x$. As an example
of that we implemented the finite difference method in Matlab for the
rectangular grid where $t_0=0, t_N=1, x_J=2$, $\Delta t = \Delta x =
0.01$, $\lambda=0.01$ and $U$ is an inverse Gaussian process with
parameters $\delta= 15$, $\gamma= 1$, and the kernel function \eqref
{modBj} has parameters $a=b=1,\alpha=0.01$. See Figure~\ref{figIGGB}
for a plot of the relative error between the boundary of the finite
difference method and the numerical integration method, and the field
obtained. For reference, we simulated the same rectangular grid by
means of numerical integration for each fixed $x$. Using the tic, toc
Matlab function, we measured the efficiency of the respective methods
in terms of speed. Unsurprisingly the finite difference method was
faster, using 0.0731 sec, whereas the numerical integration method used
0.3536 sec (the experiments were performed on a standard laptop computer).

%
%

Finally, we remark that it is also easy to estimate the error from
estimating the volatility as follows: For a constant $C = \mathbb{E}[U^2(1)]
\geq0$ and $r > s$ it holds that
\begin{eqnarray*}
\mathbb{E}\bigl[\bigl|\sigma(s) - \sigma(r)\bigr|^2\bigr] &=& \mathbb{E}
\bigl[Z(s) + Z(r) - 2\bigl(Z(s)Z(r)\bigr)^{1/2} \bigr]
\\
&= & C \biggl(\int_{-\infty}^{s}\mathrm{e}^{-\lambda(s-u)}\,\mathrm{d}u +
\int_{-\infty}^r \mathrm{e}^{-\lambda(r-u)}\,\mathrm{d}u \biggr)
\\
&&{}- 2\mathbb{E} \biggl[ \biggl(\int_{-\infty}^{s}\mathrm{e}^{-\lambda
(s-u)}\,\mathrm{d}U(u)
\int_{-\infty}^r \mathrm{e}^{-\lambda(r-u)}\,\mathrm{d}U(u)
\biggr)^{1/2} \biggr],
\end{eqnarray*}
\noindent and by non-negativity of the stochastic integral driven by a
subordinator it holds that
\begin{eqnarray*}
&& \mathbb{E}  \biggl[ \biggl(\int_{-\infty}^{s}\mathrm{e}^{-\lambda
(s-u)}\,\mathrm{d}U(u)
\int_{-\infty}^r \mathrm{e}^{-\lambda(r-u)}\,\mathrm{d}U(u)
\biggr)^{1/2} \biggr]
\\
&&\quad= \mathrm{e}^{-\lambda(r-s)/2}\mathbb{E} \biggl[ \biggl(\int_{-\infty
}^{s}\mathrm{e}^{-\lambda
(s-u)}\,\mathrm{d}U(u)
\biggl(\int_{-\infty}^{s}\mathrm{e}^{-\lambda(s-u)}\,\mathrm{d}U(u) + \int
_{s}^r\mathrm{e}^{-\lambda(s-u)}\,\mathrm{d}U(u) \biggr)
\biggr)^{1/2} \biggr]
\\
&&\quad\geq \mathrm{e}^{-\lambda(r-s)/2}\mathbb{E} \biggl[\int_{-\infty
}^{s}\mathrm{e}^{-\lambda
(s-u)}\,\mathrm{d}U(u)
\biggr].
\end{eqnarray*}
So for $r > s$, we may conclude that
\begin{eqnarray*}
\mathbb{E}\bigl[\bigl|\sigma(s) - \sigma(r)\bigr|^2\bigr] &\leq&  C \biggl(\int
_{-\infty}^r \mathrm{e}^{-\lambda(r-u)}\,\mathrm{d}u - \bigl(2\mathrm{e}^{-\lambda
(r-s)/2}
- 1\bigr)\int_{-\infty}^{s}\mathrm{e}^{-\lambda(s-u)}\,\mathrm{d}u \biggr)
\\
&= & \frac{2C}{\lambda} \bigl(1 - \mathrm{e}^{-\lambda(r-s)/2} \bigr),
\end{eqnarray*}
and thus by taking supremum we conclude that
\[
\sup_{|s-r|< \Delta t}\mathbb{E}\bigl[\bigl|\sigma(s) - \sigma(r)\bigr|^2
\bigr] \leq \frac
{2C}{\lambda} \bigl(1 - \mathrm{e}^{-\lambda\Delta t/2} \bigr).
\]

%

Having benchmarked our method of obtaining space time fields against
the more straightforward approach of numerical integration, we would
like to point out that our method has a variety of potential
applications. One might for example consider the problem of simulating
fractional Brownian motion (see, e.g., Biagini \textit{et al.} \cite{BHOZ}).
Recall that for a given Hurst parameter $H \in(0,1)$ fractional
Brownian motion can be written as
\[
B^H(t) = \frac{1}{\Gamma(H+1/2)} \biggl(\int_{-\infty}^t
(t-s)^{H-1/2}\,\mathrm{d}B(s) - \int_{-\infty}^0(-s)^{H-1/2}\,\mathrm{d}B(s)
\biggr),
\]
for $t \in\mathbb{R}$. Now notice that the kernel function $g(u) = u^{H-1/2}$
is not Lipschitz at the origin. Thus, we can not apply our convergence
result \ref{prop:Conv} directly. However, we may for a given $\varepsilon
> 0$ define an approximative kernel function
\[
h_\varepsilon(u) = \cases{
g(u), & \quad\mbox{if $u \geq \varepsilon$},\vspace*{3pt}\cr
g(\varepsilon), &\quad\mbox{if $u \in[0,\varepsilon]$,}}
\]
and employ Lemma~\ref{lem:InitCond} find that
\[
\Vert g - h_\varepsilon\Vert_{L^2(\mathbb{R}_+)}^2 \leq(2 + 1/H)
\varepsilon^{2H}.
\]
So unsurprisingly this estimate is better for $H$ closer to one than
the origin. Hence we may again employ Lemma~\ref{lem:InitCond} together
with Proposition~\ref{prop:Conv} to control simulation errors when
employing the finite difference scheme with the kernel function
$h_\varepsilon$ to simulate a trajectory of fractional Brownian motion for
a given Hurst parameter $H$.

\section{Conclusion}
\label{con}
We have defined, and analysed, a finite difference method for
simulating mild solutions of a particular HSPDE. Further we have
described how VMV processes may be viewed as mild solutions of these
particular HSPDEs, and thus obtained an alternative to numerical
integration for simulating VMV processes. Finally, we have seen in
experiments that our finite difference method is more time efficient
than numerical integration for simulating a space time random field LSS
process driven by non-exponential kernel functions. Our examples also
include the simulation of fractional Brownian/L\'evy random fields. We
remark that the finite difference scheme may also be applied for the
simulation of forward rates in the Musiela parametrisation of the
Heath--Jarrow--Morton modelling approach in fixed-income markets (see
\cite{HJM}). In future studies, we will extend our HSPDE approach to
the simulation of so-called ambit fields (see \cite{BN1}).

\section*{Acknowledgements}
We are grateful to Ole E. Barndorff-Nielsen and Almut Veraart for their
valuable suggestions, and for fruitful criticism from an anonymous
referee. Financial support from the Norwegian Research Council of the
eVita project 205328 ``Energy Markets: modeling, optimization and
simulation'' (Emmos) is greatly acknowledged. Heidar Eyjolfsson moreover
acknowledges funding from Finansmarkedsfondet.


\printhistory
\end{document}